\overfullrule=0pt
\centerline {\bf An alternative theorem for gradient systems}\par
\bigskip
\bigskip
\centerline {BIAGIO RICCERI}\par
\bigskip
\bigskip
{\bf Abstract:} In this paper, given two Banach spaces $X, Y$ and a $C^1$ functional
$\Phi:X\times Y\to {\bf R}$, under general assumptions, we show that either $\Phi$ has
a saddle-point in $X\times Y$ or, for each convex and dense set $S\subseteq Y$, there is some $\tilde y\in S$ such that
$\Phi(\cdot,\tilde y)$ has at least three critical points in $X$, two of which are global minima. Also, an application to
non-cooperative elliptic systems is presented.
\bigskip
\bigskip
\bigskip
\bigskip
The present paper is part of the extensive program of studying consequences and applications of certain general minimax theorems
([9], [10], [12]-[25]) which cannot be directly deduced by the classical Fan-Sion theorem ([5], [26]).\par
\smallskip
Here, we are interested in gradient systems. Precisely, given two Banach spaces $X, Y$ and a $C^1$ functional
$\Phi:X\times Y\to {\bf R}$, we are interested in the existence of critical points for $\Phi$, that is in the solvability
of the system
$$\cases {\Phi'_x(x,y)=0 \cr & \cr \Phi'_y(x,y)=0\ ,\cr}$$
where $\Phi'_x$ (resp. $\Phi'_y$) is the derivative of $\Phi$ with respect to $x$ (resp. $y$).\par
\smallskip
Let $I:X\to {\bf R}$. As usual, $I$ is said to be coercive if $\lim_{\|x\|\to +\infty}I(x)=+\infty$. $I$ is said to be quasi-concave (resp. quasi-convex)
if the set $I^{-1}([r,+\infty[)$ (resp. $I^{-1}(]-\infty,r])$) is convex for all $r\in {\bf R}$.
When $I$ is $C^1$, it is said to
satisfy the Palais-Smale condition if each sequence $\{x_n\}$ in $X$ such that $\sup_{n\in {\bf N}}|I(x_n)|<+\infty$ and
$\lim_{n\to \infty}\|I'(x_n)\|_{X^*}=0$ admits a strongly convergent subsequence.\par
\smallskip
Here is  our main abstract theorem:\par
\medskip
THEOREM 1. - {\it Let $X, Y$ be two real reflexive Banach spaces
 and let $\Phi:X\times Y\to {\bf R}$ be a $C^1$ functional satisfying the following
conditions:\par
\noindent
$(a)$\hskip 5pt 
the functional $\Phi(x,\cdot)$ is quasi-concave for all $x\in X$ and the functional $-\Phi(x_0,\cdot)$ is coercive for some $x_0\in X$;\par
\noindent
$(b)$\hskip 5pt  there exists a convex set $S\subseteq Y$ dense in $Y$, such that, for each $y\in S$, the functional
$\Phi(\cdot,y)$ is weakly lower semicontinuous, coercive and satisfies the Palais-Smale condition\ .\par
Then, either the system
$$\cases {\Phi'_x(x,y)=0 \cr & \cr \Phi'_y(x,y)=0\cr}$$
has a solution $(x^*,y^*)$ such that
$$\Phi(x^*,y^*)=\inf_{x\in X}\Phi(x,y^*)=\sup_{y\in Y}\Phi(x^*,y)\ ,$$
or, for every convex set $T\subseteq S$ dense in $Y$, there exists $\tilde y\in T$ such that
equation
$$\Phi'_x(x,\tilde y)=0$$ has at least three solutions, two of which are global minima in $X$ of the functional $\Phi(\cdot,\tilde y)$.}\par
\smallskip
PROOF. Assume that there is no solution $(x^*,y^*)$ of the system
$$\cases {\Phi'_x(x,y)=0 \cr & \cr \Phi'_y(x,y)=0\cr}$$
such that
$$\Phi(x^*,y^*)=\inf_{x\in X}\Phi(x,y^*)=\sup_{y\in Y}\Phi(x^*,y)\ .$$
We consider both $X, Y$ endowed with the weak topology. Notice that, by $(a)$,
$\Phi(x,\cdot)$ is weakly upper semicontinuous in $Y$ for all $x\in X$ and weakly $\sup$-compact for $x=x_0$. As a consequence, the functional
$y\to \inf_{x\in X}\Phi(x,y)$ is weakly $\sup$-compact and so it attains its supremum.
Likewise, by $(b)$,
$\Phi(\cdot,y)$ is weakly $\inf$-compact for all $y\in S$. By continuity and density, we have
$$\sup_{y\in Y}\Phi(x,y)=\sup_{y\in S}\Phi(x,y)\eqno{(1)}$$
for all $x\in X$. As a consequence, the functional $x\to \sup_{y\in Y}\Phi(x,y)$ is weakly $\inf$-compact and so it attains its infimum.
Therefore, the occurrence of the equality
$$\sup_Y\inf_X\Phi=\inf_X\sup_Y\Phi$$
is equivalent to the existence of a point $(\hat x,\hat y)\in X\times Y$ such that
$$\sup_{y\in Y}\Phi(\hat x,y)=\Phi(\hat x,\hat y)=\inf_{x\in X}\Phi(x,\hat y)\ .$$
But, for what we are assuming, no such a point can exist and hence we have
$$\sup_Y\inf_X\Phi<\inf_X\sup_Y\Phi\ .\eqno{(2)}$$
So, in view of $(1)$ and $(2)$, we also have
$$\sup_S\inf_X\Phi<\inf_X\sup_S\Phi\ .$$
At this point, we are allowed to apply Theorem 1.1 of [20]. Therefore, there exists $\tilde y\in S$ such that the functional $\Phi(\cdot,\tilde y)$ has at
least two global minima in $X$ and so, thanks to Corollary 1 of [8], the same functional has at least three critical points.\hfill $\bigtriangleup$\par
\medskip
The next result is a consequence of Theorem 1.\par
\medskip
THEOREM 2. - {\it Let $X, Y$ be two real Hilbert spaces and let $J:X\times Y\to {\bf R}$ be a $C^1$ functional satisfying the following conditions:\par
\noindent
$(a_1)$\hskip 5pt  the functional $y\to {{1}\over {2}}\|y\|_Y^2+J(x,y)$ is quasi-convex for all $x\in X$ and coercive for some $x\in X$\ ;\par
\noindent
$(b_1)$\hskip 5pt there exists a convex set $S\subseteq Y$ dense in $Y$ such that, for each $y\in S$, the operator $J'_x(\cdot,y)$ is compact and
$$\limsup_{\|x\|_X\to +\infty}{{J(x,y)}\over {\|x\|_X^2}}<{{1}\over {2}}\ ;\eqno{(3)}$$
Then,  either the system
$$\cases {x=J'_x(x,y) \cr & \cr y=-J'_y(x,y)\cr}$$
has a solution $(x^*,y^*)$ such that
$${{1}\over {2}}(\|x^*\|_X^2-\|y^*\|_Y^2)-J(x^*,y^*)=
\inf_{x\in X}\left ({{1}\over {2}}(\|x\|_X^2-\|y^*\|_Y^2)-J(x,y^*)\right )=
\sup_{y\in Y}\left ({{1}\over {2}}(\|x^*\|_X^2-\|y\|_Y^2)-J(x^*,y)\right )\ ,$$
or, for every convex set $T\subseteq S$ dense in $Y$,
there exists $\tilde y\in T$ such that
the equation
$$x=J'_x(x,\tilde y)$$
has at least three solutions, two of which are global minima in $X$ of the functional $x\to {{1}\over {2}}\|x\|_X^2-J(x,\tilde y)$.}\par
\smallskip
PROOF.   Consider the function $\Phi:X\times Y\to {\bf R}$ defined by
$$\Phi(x,y)={{1}\over {2}}(\|x\|_X^2-\|y\|_Y^2)-J(x,y)$$
for all $(x,y)\in X\times Y$. Clearly, $\Phi$ is $C^1$ and one has
$$\Phi'_x(x,y)=x-J'_x(x,y)\ ,$$
$$\Phi'_y(x,y)=-y-J'_y(x,y)$$
for all $(x,y)\in X\times Y$. We want to apply Theorem 1 such a $\Phi$. Of course, $\Phi$ satisfies $(a)$ in view of $(a_1)$. Concerning $(b)$,
notice that, for each $y\in S$, the functional $J(\cdot,y)$ is sequentially weakly continuous since $J'_x(\cdot,y)$ is compact ([27], Corollary 41.9).
Moreover, from $(3)$ it immediately follows that $\Phi(\cdot,y)$ is coercive and so, by the Eberlein-Smulyan theorem, it is weakly lower semicontinuous.
Finally, $\Phi(\cdot,y)$ satisfies the Palais-Smale condition in view of Example 38.25 of [27]. Now, the conclusion follows directly from Theorem 1.
\hfill $\bigtriangleup$\par
\medskip
We now present an application of Theorem 2  to non-cooperative elliptic systems.\par
\smallskip
In what follows, $\Omega\subset {\bf R}^n$ ($n\geq 2$) is a bounded smooth domain. We consider $H^1_0(\Omega)$ equipped with the
scalar product
$$\langle u,v\rangle=\int_{\Omega}\nabla u(x)\nabla v(x)dx\ .$$
\smallskip
We denote by ${\cal A}$ the class of all functions $H:\Omega\times {\bf R}^2\to {\bf R}$, with $H(x,0,0)=0$ for all $x\in \Omega$, which are measurable in $\Omega$,
$C^1$ in ${\bf R}^2$ and satisfy 
$$\sup_{(x,u,v)\in \Omega\times {\bf R}^2}{{|H_u(x,u,v)|+|H_v(x,u,v)|}\over {1+|u|^p+|v|^q}}<+\infty$$
where $p, q>0$, with $p<{{n+2}\over {n-2}}$ and $q\leq {{n+2}\over {n-2}}$ when $n>2$. 
\medskip
Given $H\in {\cal A}$, we are interested in the problem\par
$$\cases {-\Delta u=H_u(x,u,v) & in $\Omega$
\cr & \cr
-\Delta v=-H_v(x,u,v) & in $\Omega$
\cr & \cr
u=v=0 & on $\partial\Omega$\ ,\cr} \eqno{(P_H)}$$
$H_u$ (resp. $H_v$) denoting the derivative of $H$ with respect to $u$ (resp. $v$).\par
\smallskip
As usual, a weak solution of $(P_H)$ is any $(u,v)\in H^1_0(\Omega)\times H^1_0(\Omega)$ such that
$$\int_{\Omega}\nabla u(x)\nabla\varphi(x)dx=\int_{\Omega}H_u(x,u(x),v(x))\varphi(x)dx\ ,$$
$$\int_{\Omega}\nabla v(x)\nabla\psi(x)dx=-\int_{\Omega}H_v(x,u(x),v(x))\psi(x)dx$$
for all $\varphi, \psi\in H^1_0(\Omega)$.\par
\smallskip
Define the functional $I_H:H^1_0(\Omega)\times H^1_0(\Omega)\to {\bf R}$ by
$$I_H(u,v)={{1}\over {2}}\left ( \int_{\Omega}|\nabla u(x)|^2dx-\int_{\Omega}|\nabla v(x)|^2dx\right )-
\int_{\Omega}H(x,u(x),v(x))dx$$
for all $(u,v)\in H^1_0(\Omega)\times H^1_0(\Omega)$.\par
\smallskip
Since $H\in {\cal A}$, the functional $I_H$ is $C^1$ in $H^1_0(\Omega)\times H^1_0(\Omega)$ and its critical
points are precisely the weak solutions of $(P_H)$. 
\smallskip
Also, we denote by $\lambda_1$ the first eigenvalue of the Dirichlet problem
$$\cases{-\Delta u=\lambda u & in $\Omega$ \cr & \cr u=0 & on $\partial\Omega$\ .\cr}$$
Our result on $(P_H)$ is a follows:\par
\medskip
THEOREM 3. - {\it Let $H\in {\cal A}$ be such that 
$$\limsup_{|u|\to +\infty}{{\sup_{x\in \Omega}\sup_{|v|\leq r}H(x,u,v)}\over {u^2}}\leq 0 \eqno{(4)}$$
for all $r>0$,
and
$$\hbox {\rm meas}\left (\left\{x\in \Omega : \sup_{u\in {\bf R}}H(x,u,0)>0\right\}\right)>0\ . \eqno{(5)}$$
Moreover, assume that either $H(x,u,\cdot)$ is convex for all $(x,u)\in \Omega\times {\bf R}$, or
$$L:=\sup_{(v,\omega)\in {\bf R}^2, v\neq \omega}{{\sup_{(x,u)\in \Omega\times {\bf R}}|H_v(x,u,v)-H_v(x,u,\omega)|}\over {|v-\omega|}}<+\infty\ .\eqno{(6)}$$
Set
$$\lambda^*={{1}\over {2}}\inf\left\{{{\int_{\Omega}|\nabla w(x)|^2dx}\over {\int_{\Omega}H(x,w(x),0)dx}} : w\in H^1_0(\Omega), \int_{\Omega}H(x,w(x),0)dx>0\right\}$$
and assume that $\lambda^*<{{\lambda_1}\over {L}}$ when $(6)$ holds.\par
Then, for each $\lambda>\lambda^*$, with $\lambda< {{\lambda_1}\over {L}}$ when $(6)$ holds, either the problem
$$\cases {-\Delta u=\lambda H_u(x,u,v) & in $\Omega$
\cr & \cr
-\Delta v=-\lambda H_v(x,u,v) & in $\Omega$
\cr & \cr
u=v=0 & on $\partial\Omega$\cr} $$
has a non-zero weak solution belonging to $L^{\infty}(\Omega)\times L^{\infty}(\Omega)$, or, for each convex set $T\subseteq H^1_0(\Omega)\cap L^{\infty}(\Omega)$ dense in $H^1_0(\Omega)$, there exists
$\tilde v\in T$ such that the problem
$$\cases {-\Delta u=\lambda H_u(x,u,\tilde v(x)) & in $\Omega$
\cr & \cr
u=0 & on $\partial\Omega$\cr}$$
has at least three weak solutions, two of which are global minima in $H^1_0(\Omega)$ of the functional $I_{\lambda H}(\cdot,\tilde v)$.}\par
\smallskip
PROOF. Define the functional $J:H^1_0(\Omega)\times H^1_0(\Omega)\to {\bf R}$ by
$$J(u,v)=\int_{\Omega}H(x,u(x),v(x))dx$$
for all $(u,v)\in H^1_0(\Omega)\times H^1_0(\Omega)$. Notice that $(5)$ implies $\sup_{u\in H^1_0(\Omega)}J(u,0)>0$ ([16], pp. 135-136). Consequently, $\lambda^*<+\infty$.
Fix $\lambda>\lambda^*$, with $\lambda< {{\lambda_1}\over {L}}$ when $(6)$ holds.
 We want to apply Theorem 2 to $\lambda J$. Concerning $(a_1)$, notice that, for each $u\in H^1_0(\Omega)$,
the functional $v\to {{1}\over {2}}\int_{\Omega}|\nabla v(x)|^2dx+\lambda J(u,v)$ is strictly convex and coercive. This is clear when $H(x,\xi,\cdot)$ is convex for all
$(x,\xi)\in {\bf R}^2$. When $(6)$ holds, the operator $J'_v(u,\cdot)$ turns out to be Lipschitzian in $H^1_0(\Omega)$ with Lipschitz constant
${{L}\over {\lambda_1}}$  ([11], p. 165). So, the operator $v\to v-\lambda J'_v(u,v)$ is uniformly monotone and
then the claim follows from a classical result ([27], pp. 247-249). Concerning $(b_1)$, fix $v\in H^1_0(\Omega)\cap L^{\infty}(\Omega)$. Notice that 
$J'_u(\cdot,v)$ is compact due to restriction on $p$ (recall that $H\in {\cal A}$). Moreover, in view of $(4)$, for each $\epsilon>0$, there exists $\delta>0$ 
such that
$$H(x,t,s)\leq \epsilon t^2$$
for all $x\in \Omega$, $s\in \left [-\|v\|_{L^{\infty}(\Omega)},\|v\|_{L^{\infty}(\Omega)}\right ]$ and $t\in {\bf R}\setminus [-\delta,\delta]$. But $H$ is bounded
on each bounded subset of $\Omega\times {\bf R}^2$, and so, for a suitable constant $c>0$, we have
$$H(x,t,s)\leq \epsilon t^2+c \eqno{(7)}$$
for all $(x,t,s)\in \Omega\times {\bf R}\times \left [-\|v\|_{L^{\infty}(\Omega)},\|v\|_{L^{\infty}(\Omega)}\right ]$. Of course, from $(7)$ it follows that
$$\limsup_{\|u\|\to +\infty}{{J(u,v)}\over {\|u\|^2}}\leq \epsilon$$
and so 
$$\limsup_{\|u\|\to +\infty}{{J(u,v)}\over {\|u\|^2}}\leq 0$$
since $\epsilon>0$ is arbitrary. Hence, $\lambda J$ satisfies $(3)$. Now suppose that there exists a convex set 
$T\subseteq H^1_0(\Omega)\cap L^{\infty}(\Omega)$ dense in $H^1_0(\Omega)$ such that, for each $v\in T$, the problem
$$\cases {-\Delta u=\lambda H_u(x,u,v(x)) & in $\Omega$
\cr & \cr
u=0 & on $\partial\Omega$\cr}$$
has at most two weak solutions. Then, Theorem 2 ensures the existence of a weak solution $(u^*,v^*)$ of the problem
$$\cases {-\Delta u=\lambda H_u(x,u,v) & in $\Omega$
\cr & \cr
-\Delta v=-\lambda H_v(x,u,v) & in $\Omega$
\cr & \cr
u=v=0 & on $\partial\Omega$\cr} $$
such that 
$$I_{\lambda H}(u^*,v^*)=\inf_{u\in H^1_0(\Omega)}I_{\lambda H}(u,v^*)=\sup_{v\in H^1_0(\Omega)}I_{\lambda H}(u^*,v)\ .\eqno{(8)}$$
From $(8)$, in view of Theorem 1 of [3] (see Remark 5, p. 1631), it follows that $u^*, v^*\in L^{\infty}(\Omega)$. We show that $(u^*,v^*)\neq
(0,0)$. If $v^*\neq 0$, we are done. So, assume $v^*=0$. Since $\lambda>\lambda^*$, we have
$$\inf_{u\in H^1_0(\Omega)}\left ({{1}\over {2}}\int_{\Omega}|\nabla u(x)|^2dx-\lambda\int_{\Omega}H(x,u(x),0)dx\right )<0\ .\eqno{(9)}$$
But then, since $\int_{\Omega}H(x,0,0)dx=0$, from $(9)$ and the first equality in $(8)$, it follows that $u^*\neq 0$, and the proof is complete.
\hfill $\bigtriangleup$\par
\medskip
For previous results on problem $(P_H)$ (markedly different from Theorem 3) we refer to [1], [4], [6], [7].\par
\smallskip
A joint application of Theorem 3 with the main result in [2] gives the following:\par
\medskip
THEOREM 4. - {\it Let $H\in {\cal A}$ satisfy the assumptions of Theorem 3. Moreover, suppose that
 $\inf_{\Omega\times {\bf R}^2} H_u\geq 0$  and that,
 for each $(x,v)\in \Omega\times {\bf R}$, the function $u\to {{H_u(x,u,v)}\over {u}}$ is strictly decreasing in $]0,+\infty[$.\par
Then, for every $\lambda>\lambda^*$, with $\lambda< {{\lambda_1}\over {L}}$ when $(6)$ holds, the problem
$$\cases {-\Delta u=\lambda H_u(x,u,v) & in $\Omega$
\cr & \cr
-\Delta v=-\lambda H_v(x,u,v) & in $\Omega$
\cr & \cr
u=v=0 & on $\partial\Omega$\cr} $$
has a non-zero weak solution belonging to $L^{\infty}(\Omega)\times L^{\infty}(\Omega)$.}\par
\smallskip 
PROOF. Fix $\lambda>\lambda^*$, with $\lambda< {{\lambda_1}\over {L}}$ when $(6)$ holds. Fix also $v\in C^{\infty}_0(\Omega)$. Since $\inf_{\Omega\times {\bf R}^2} H_u\geq 0$,
the bounded weak solutions of the problem
$$\cases {-\Delta u=\lambda H_u(x,u,v(x)) & in $\Omega$
\cr & \cr
u=0 & on $\partial\Omega$\cr}$$
are continuous and non-negative in $\overline {\Omega}$. As a consequence, in view of Theorem 1 of [2], the problem
$$\cases {-\Delta u=\lambda H_u(x,u,v(x)) & in $\Omega$
\cr & \cr
u=0 & on $\partial\Omega$\cr}$$
has at most one non-zero bounded weak solution. Now, the conclusion follows directly from Theorem 3.\hfill $\bigtriangleup$\par
\medskip
Finally, notice the following corollary of Theorem 4:\par
\medskip
THEOREM 5. - {\it Let $F, G : {\bf R}\to {\bf R}$ be two $C^1$ functions, with $FG-F(0)G(0)\in {\cal A}$, satisfying the following conditions:\par
\noindent
$(a_2)$\hskip 5pt $F$ is non-negative, increasing, $\lim_{u\to +\infty}{{F(u)}\over {u^2}}=0$ and the function $u\to {{F'(u)}\over {u}}$ is strictly decreasing in
$]0,+\infty[$\ ;\par
\noindent
$(b_2)$\hskip 5pt $G$ is positive and convex.\par
Finally, let $\alpha\in L^{\infty}(\Omega)$, with $\alpha>0$. Set
$$\lambda_{\alpha}^*={{1}\over {2G(0)}}\inf\left\{{{\int_{\Omega}|\nabla w(x)|^2dx}\over {\int_{\Omega}\alpha(x)(F(w(x))-F(0))dx}} : w\in H^1_0(\Omega), \int_{\Omega}\alpha(x)(F(w(x))-F(0))dx>0\right\}\ .$$
Then, for every $\lambda>\lambda_{\alpha}^*$, the problem
$$\cases {-\Delta u=\lambda\alpha(x) G(v(x))F'(u) & in $\Omega$
\cr & \cr
-\Delta v=-\lambda\alpha(x) F(u(x))G'(v) & in $\Omega$
\cr & \cr
u=v=0 & on $\partial\Omega$\cr} $$
has a non-zero weak solution belonging to $L^{\infty}(\Omega)\times L^{\infty}(\Omega)$.}\par
\smallskip
PROOF. Apply Theorem 4 to the function $H:\Omega\times {\bf R}^2\to {\bf R}$ defined by
$$H(x,u,v)=\alpha(x)(F(u)G(v)-F(0)G(0))$$
for all $(x,u,v)\in \Omega\times {\bf R}^2$. Checking that $H$ satisfies the assumptions of Theorem 4 is an easy task.\hfill $\bigtriangleup$.

\bigskip
\bigskip
{\bf Acknowledgement.} The author has been supported by the Gruppo Nazionale per l'Analisi Matematica, la Probabilit\`a e 
le loro Applicazioni (GNAMPA) of the Istituto Nazionale di Alta Matematica (INdAM) and by the Universit\`a degli Studi di Catania, ``Piano della Ricerca 2016/2018 Linea di intervento 2". \par
\bigskip
\bigskip
\bigskip
\bigskip
\centerline {\bf References}\par
\bigskip
\bigskip
\noindent
[1]\hskip 5pt J. C. BATKAM and F. COLIN, {\it The effects of concave and convex nonlinearities in some noncooperative elliptic systems}, 
Ann. Mat. Pura Appl., {\bf 193} (2014), 1565-1576.\par
\smallskip
\noindent
[2]\hskip 5pt H. BREZIS and L. OSWALD, {\it Remarks on sublinear elliptic equation}, Nonlinear Anal., {\bf 10} (1986), 55-64.\par
\smallskip
\noindent
[3]\hskip 5pt A. CIANCHI, {\it Boundedness of solutions to variational problems under general growth conditions},  Comm. Partial Differential Equations,
 {\bf 22} (1997), 1629-1646.
\smallskip
\noindent
[4]\hskip 5pt D. G. DE FIGUEIREDO and Y. H. DING, {\it  Strongly indefinite functionals and multiple solutions of elliptic systems}, Trans. Amer. Math. Soc., {\bf 355} (2003), 2973-2989.\par
\smallskip
\noindent
[5]\hskip 5pt K. FAN, {\it Fixed-point and minimax theorems in locally convex topological
linear spaces}, Proc. Nat. Acad. Sci. U.S.A., {\bf 38} (1952), 121-126.\par
\smallskip
\noindent
[6]\hskip 5pt Y. GUO, {\it Nontrivial solutions for resonant noncooperative elliptic systems}, Comm. Pure Appl. Math., {\bf 53} (2000), 1335-1349.\par
\smallskip
\noindent
[7]\hskip 5pt N. HIRANO, {\it Infinitely many solutions for non-cooperative elliptic systems}, J. Math. Anal. Appl., {\bf 311} (2005), 545-566.\par
\smallskip
\noindent
[8]\hskip 5pt P. PUCCI and J. SERRIN, {\it A mountain pass
theorem}, J. Differential Equations, {\bf 60} (1985), 142-149.\par
\smallskip
\noindent
[9]\hskip 5pt B. RICCERI, {\it Some topological mini-max theorems via
an alternative principle for multifunctions}, Arch. Math. (Basel),
{\bf 60} (1993), 367-377.\par
\smallskip
\noindent
[10]\hskip 5pt B. RICCERI, {\it On a topological minimax theorem and
its applications}, in ``Minimax theory and applications'', B. Ricceri
and S. Simons eds., 191-216, Kluwer Academic Publishers, 1998.\par
\smallskip
\noindent
[11]\hskip 5pt B. RICCERI, {\it Sublevel sets and global minima of coercive
functionals and local minima of their perturbations}, J. Nonlinear Convex
Anal., {\bf 5} (2004), 157-168.\par
\smallskip
\noindent
[12]\hskip 5pt B. RICCERI, {\it Minimax theorems for functions involving a
real variable and applications}, Fixed Point Theory, {\bf 9} (2008),
275-291.\par
\smallskip
\noindent
[13]\hskip 5pt B. RICCERI, {\it Well-posedness of constrained minimization
problems via saddle-points}, J. Global Optim., {\bf 40} (2008),
389-397.\par
\smallskip
\noindent
[14]\hskip 5pt B. RICCERI, {\it Multiplicity of global minima for
parametrized functions}, Rend. Lincei Mat. Appl., {\bf 21} (2010),
47-57.\par
\smallskip
\noindent
[15]\hskip 5pt B. RICCERI, {\it A strict minimax inequality criterion and some of its consequences}, Positivity, {\bf 16} (2012), 455-470.\par
\smallskip
\noindent
[16]\hskip 5pt B. RICCERI, {\it Energy functionals of Kirchhoff-type problems having multiple global minima}, Nonlinear Anal., {\bf 115} (2015),  
130-136.\par
\smallskip
\noindent
[17]\hskip 5pt B. RICCERI, {\it A minimax theorem in infinite-dimensional topological vector spaces}, Linear Nonlinear Anal., {\bf 2}
(2016), 47-52.\par
\smallskip
\noindent
[18]\hskip 5pt B. RICCERI, {\it Miscellaneous applications of certain minimax theorems I}, Proc. Dynam. Systems Appl., {\bf 7} (2016),
198-202.\par
\smallskip
\noindent
[19]\hskip 5pt B. RICCERI, {\it On the infimum of certain functionals}, in ``Essays in Mathematics and its Applications -
In Honor of Vladimir Arnold", Th. M. Rassias and P. M. Pardalos eds., 361-367, Springer, 2016.\par
\smallskip
\noindent
[20]\hskip 5pt B. RICCERI, {\it On a minimax theorem: an improvement, a new proof and an overview of its applications},
Minimax Theory Appl., {\bf 2} (2017), 99-152.\par
\smallskip
\noindent
[21]\hskip 5pt B. RICCERI, {\it Miscellaneous applications of certain minimax theorems II}, Acta Math. Vietnam.,  {\bf 45} (2020), 515-524.\par
\smallskip
\noindent
[22]\hskip 5pt B. RICCERI, {\it Minimax theorems in a fully non-convex setting}, J. Nonlinear Var. Anal., {\bf 3} (2019), 45-52.\par
\smallskip
\noindent
[23]\hskip 5pt B. RICCERI, {\it Applying twice a minimax theorem}, J. Nonlinear Convex Anal., {\bf 20} (2019), 1987-1993.\par
\smallskip
\noindent
[24]\hskip 5pt B. RICCERI, {\it Another multiplicity result for the periodic solutions of certain systems}, Linear Nonlinear Anal., {\bf 5} (2019),
371-378.\par
\smallskip
\noindent
[25]\hskip 5pt B. RICCERI, {\it A remark on variational inequalities in small balls}, J. Nonlinear Var. Anal., {\bf 4} (2020), 21-26.\par
\smallskip
\noindent
[26]\hskip 5pt M. SION, {\it On general minimax theorems}, Pacific J. Math., {\bf 8} (1958), 171-176.\par
\smallskip
\noindent
[27]\hskip 5pt E. ZEIDLER, {\it Nonlinear functional analysis and its
applications}, vol. III, Springer-Verlag, 1985.\par
\bigskip
\bigskip
\bigskip
\bigskip
Department of Mathematics and Informatics\par
University of Catania\par
Viale A. Doria 6\par
95125 Catania, Italy\par
{\it e-mail address}: ricceri@dmi.unict.it

\bye